\documentclass[12pt]{amsart}
\usepackage{amsmath}
\usepackage{amsfonts,amssymb,amsthm, txfonts, pxfonts, boxedminipage}

\newtheorem{theorem}{Theorem}
\newtheorem{lemma}{Lemma}
\newtheorem{corollary}{Corollary}
\theoremstyle{remark}

\begin{document}

\title[Delaunay triangulations of extra-large metrics]{Extra-large metrics}
\author{Igor Rivin}
\begin{abstract}
We show that every two dimensional spherical cone metric with all cone
angles greater than $2\pi$ and the lengths of all closed geodesics
greater than $2\pi$ admits a triangulation whose $0$-skeleton is
precisely the set of cone points -- this is, in fact, the Delaunay
triangulation of the set of cone points.
\end{abstract}

\date{\today}
\keywords{triangulation, extra-large, hyperbolic polyhedra}
\subjclass{52B11, 52B10, 57M50}

\maketitle
\section{Introduction}

In this note we study \emph{extra large} spherical cone manifolds in
dimension $2$ (though many of our resultsz and techniques extend to
higher dimensions.

A $2$-dimensional \emph{spherical cone manifold} is a metric space
where all but finitely many points has a neighborhood isometric to a
neighborhood of a point on the round sphere $\mathbb{S}^2.$ The
exceptional points (\emph{cone points}) have neighboroods isometric to
a spherical cone, the angle of which is the \emph{cone angle} at that
point. 

If $M$ is a cone manifold, we define a \emph{geodesic} to be a locally
length minimizing curve on $M.$ It is easy to see that such a curve is
locally a great circle, except at the cone points. There, the geodesic
must have the property that it subtends an angle no smaller than $\pi$
on either side. Consequently, no geodesics can \emph{pass through}
cone points, where the cone angles are smaller than $2\pi$ (such cone
points are known as \emph{positively curved} cone points, since the
curvature of a cone point is defined as $2\pi$ less the cone angle at
the point). 

We say that a spherical cone manifold is \emph{extra large} if

\begin{enumerate}
\item All the cone points are negatively curved.
\item All closed geodesics are longer than $2\pi.$
\end{enumerate}
Such spaces are of considerable importance in geometry in general (due
to work of A.D.Aleksandrov,  M. Gromov, and then R. Charney and
M. Davis \cite{charneydavis}), and
in three-dimensional hyperbolic geometry in particular, due in large
part to the results of the author (\cite{thes0,thes,topology,combopt,ann2}), who showed that the polar duals of
convex compact polyhedra in $\mathbb{H}^3$ are precisely the extra
large spherical cone manifolds homeomorphic to $\mathbb{S}^2.$ In that
work, the term \emph{extra large} was not used -- it was invented by
Gromov, to describe the vertex links in negatively curved spaces.

The main objective of this paper is to show:
\begin{theorem}
\label{mainthm}
An extra large spherical cone surface admits a cell decomposition
whose $0$-skeleton is precisely the set of cone points.
\end{theorem}

For the impatient reader, we first give the recipe for constructing
the cell decomposition whose existence is postulated in Theorem
\ref{mainthm}:

First, recall that the \emph{Voronoi diagram} $\mathcal{V}_P$ of a metric space $M$
with respect to a point set $P=\{p_1, \dotsc, p_n, \dotsc\}$ is the
decomposition of $M$ into \emph{Voronoi cells}
\[
V_i = \{ x\in M \left| d(x, p_i) \leq d(x, p_j), \quad \forall j\right.\}.
\]
Clearly, 
\[
\bigcup_i V_i = M,
\]
and 
\[
\overset{\circ}{V}_i \cap \overset{\circ}{V}_j = \emptyset,\quad i \neq j.
\]
As will be shown below, each $V_i$ is a geodesic polygon, and the
Voronoi diagram $\mathcal{V}_P$  is a cell decomposition of $M..$ The
\emph{Delaunay tesselation} $\mathcal{D}_P$ of $M$ with respect to $P$
is the Poincar\'e dual of $\mathcal{V}_P:$  its edges corrspond to
pairs $p_i, p_j$ of sites whose Voronoi cells share an edge, while its
faces correspond to points of $M$ equidistant from three or more
elements of $P.$ The cells of $\mathcal{D}_P$ are convex, and so the
tesselation $\mathcal{D}_P$ can be completed to a triangulation of
$M.$

To push the program above through, we will need a number of steps.

\section{The injectivity radius}
\label{geomprelim}
We define the \emph{injectivity radius} of the space $M$ \emph{at a point
$p$} as the radius of the largest disk in the tangent space of $p$ for
which the exponential map is an embedding. The injectivity radius of
$M$ is the infimum over all points $p$ of the injectivity radii of $M$
at $p.$ In simpler terms, the injectivity radius of $M$ is the
smallest $d$ such that there exist at least two distinct curves from
$p$ to $q$ realizing $d(p, q) = d.$

Our first result is:
\begin{theorem}
\label{injthm}
A space $M$ is extra large if and only if the cone points of $M$ are
negatively curved and the injectivity radius of $M$ is greater than
$\pi.$
\end{theorem}
\begin{proof}
Let $p, q$ be the pair of points realizing the injectivity
radius. This means that there are two shortest curves $\gamma_1,
\gamma_2$ of length $L = \ell(\gamma_1) = \ell(\gamma_2) \leq \pi$
connecting $p$ to $q,$ and $L$ is the smallest with this property.

Let $\gamma = \gamma_1 \cup \gamma_2.$ If $\gamma$ is geodesic, then
$L(\gamma) \leq 2\pi,$ so we have a contradiction to extra-largeness. 
If not, suppose (without loss of generality) that $\gamma$ has a
``corner'' at $p.$ That means that on one side, the angle $\alpha$
subtended by $\gamma$ at $p$ is smaller than $\pi.$ If $\gamma_1$ and
$\gamma_2$ are both smooth, then take the bisector of $\alpha$ and
move $p$ a very small distance $\rho$ along the bisector. By elementary
spherical eometry,
\[
\dfrac{d\ell(\gamma_1)}{d\rho} = \dfrac{d \ell(\gamma_2)}{d\rho} < 0,
\]
which contradicts the minimality of $L.$
If (without loss of generality) $\gamma_1$ is \emph{not} smooth, while
$\gamma_2$ is, let $x$ be the cone point of $\gamma_1$ closest to $p.$
Note that if the angle subtended by $\gamma_1$ at $x$ on the side of
the corner at $p$ equals $\pi,$ then $x$ can be treated as a smooth
point, so the correct definition of $x$ is: the closest point of
$\gamma_1$ to $p,$ where $\gamma_1$ is not smooth on the side of the
corner (if such a point does not exist between $p$ and $q,$ then we
find ourselves back in the smooth case, which corresponds to $x =
q.$). Let $L_x = d(p, x).$

In any event, now, instead of the bisector of the angle $\alpha$ at
$p,$ we pick a direction, such that 
\[
\dfrac{d\ell(\gamma_2)}{d\rho} = \dfrac{d L_x}{d\rho} < 0.
\]
Such a direction exists by the intermediate value theorem.
The above argument adapts in the obvious way if both $\gamma_1$ and
$\gamma_2$ are singular.
\end{proof}

\section{Voronoi diagrams}

\begin{lemma}
\label{distbound}
For all $x \in V_i,$ $d(x, p_i) < \pi/2.$.
\end{lemma}

\begin{proof}
Suppose that there exists an $x$ contradicting the assertion of the
Lemma. Then the distance from $x$ to the cone locus of $M$ is at least
$\pi/2,$ and so there is a smooth hemisphere around $x.$ The boundary
of that hemisphere is a closed geodesic of length $2\pi.$
\end{proof}

\begin{corollary}
The diameter of the Voronoi cell $V_i$ is less than $\pi.$
\end{corollary}
We will need the following simple lemma from spherical geometry:

\begin{theorem}
The boundary of a Voronoi cell $V_i$ is a convex polygonal curve.
\end{theorem}
\begin{proof}
Let $x\in V_i.$ Let $r = d(x, p_i);$ we know that $r < \pi/2.$
Consider the disk $D_x(r)$ of radius $r$ around $x.$

There are the following possibilities:

Firstly, $p_i$ might be the \emph{only} cone point in $D_r(x).$ In
that case, a neighborhood of $x$ is in $V_i,$ and so $x$ is in the
interior of $V_i.$

Secondly, there msy be exactly one other point $p_j$ such that $d(,
p_j) = r.$ In that case, a small geodesic segment bisecting the angle
$p_i x p_j$ lies in $V_i \cap V_j.$

Thirdly, there can be a number of points $pi, p_{j_1}, p_{j_2},
\dotsc, p_{j_k}$ at distance $r$ from $x.$ In that case a small part
of the cone from $x$ to the Voronoi region of $p_i$ on the boundary of
$D_r(x)$ lies in $V_i$ -- note that this argument works in arbitrary
dimension. 
\end{proof}
\begin{theorem}
A Voronoi cell $V_i$ is star-shaped with respect to $p_i.$
\end{theorem}
\begin{proof}
Let $x\in V_i,$ and let $y$ be on the segment $p_i x.$ By the triangle
inequality, we see that
\[
d(y, \partial D_x(r)) + d(y, x) > r.
\]
Since for any $j\neq i$ we have that
\[
d(y, p_j) > d(y, \partial D_x(r)),
\]
the assertion of the Theorem follows.
\end{proof}

\begin{theorem}
Every Voronoi cell $V_i$ is convex.
\end{theorem}
\begin{proof}
By the preceding result, every Voronoi cell $V_i$ is a stashaped subset of a
cone of radius $\pi/2$ centered on $p_i,$ with geodesically convex
boundary. Take two points $p$ and $q$ in $V_i.$ If one of the angles
$p p_i q$ does not exceed $\pi,$ the result follows from elementary
spherical geometry. If both the angles $p p_i q$ are at least $\pi,$
the broken line $p p_i q$ is geodesic.
\end{proof}

\begin{theorem}
Let $V_i, V_j$ be two Voronoi cells, then $V_i \cap V_j$ is
connected. 
\end{theorem}
\begin{proof}
Let $p, q \in V_i \cap V_j.$ here is a shortest geodesic $\gamma_i$
from $p$ to $q$ in $V_i$ and a shortest geodesic $\gamma_j$ from $p$
to $q$ in $V_j.$ Since the diameters of $V_i$ and $V_j$ are smaller
than $\pi,$ it follows that $\gamma_i = \gamma_j.$ Thus, $\gamma =
\gamma_i = \gamma_j \subseteq V_i \cap V_j,$ hence $V_i \cap V_j$ is
path connected. In fact, the argument (together with the results
above) easily shows that $V_i \cap V_j$ is an edge of both.
\end{proof}

The above results sow that $\mathcal{V}_P = \{V_1, \dotsc, V_n,
\dotsc\}$ is a simplicial cell decomposition, and so its dual is a
cellulation of $M$ with vertices at $p_1, \dotsc, p_n, \dotsc.$ The
cells of cellulation are convex (in fact inscribed in circles; The
centers are precisely the corners of the boundaries of the cells
$V_i.$), and so the proof of Theorem \ref{mainthm} is complete.

\section{Remarks}
An identical argument with the appropriate modification of the
extra-largeness hypothesis can be used to show an analogous result for
a Riemannian surface with cone singularies. In particular, for
Euclidean and Hyperbolic cone surfaces, it is sufficient to require
the cone angles to be non-positively curved.

\bibliographystyle{plain}
\bibliography{curves,rivin,opt,matrix}

\begin{thebibliography}{1}

\bibitem{charneydavis}
Ruth Charney and Michael Davis.
\newblock The polar dual of a convex polyhedral set in hyperbolic space.
\newblock {\em Michigan Math. J.}, 42(3):479--510, 1995.

\bibitem{thes0}
Igor Rivin.
\newblock {\em On the Geometry of Convex Polyhedra in Hyperbolic 3-Space}.
\newblock PhD thesis, Princeton University, July 1986.

\bibitem{topology}
Igor Rivin.
\newblock On the geometry of ideal polyhedra in hyperbolic 3-space.
\newblock 32(1):87--92, January 1993.

\bibitem{ann2}
Igor Rivin.
\newblock Euclidean structures on simplicial surfaces and hyperbolic volume.
\newblock {\em Annals of Mathematics (ser. 2)}, 139(3):553--580, May 1994.

\bibitem{combopt}
Igor Rivin.
\newblock Combinatorial optimization in geometry.
\newblock {\em Advances in Applied Mathematics}, 31(1):242--271, 2003.
\newblock arxiv.org preprint math.GT/9907032.

\bibitem{thes}
Igor Rivin and C.D.Hodgson.
\newblock A characterization of compact convex polyhedra in hyperbolic 3-space.
\newblock {\em Inventiones Mathematicae}, pages 77--111, January 1993.
\newblock Corrigendum, vol 117, page 359.

\end{thebibliography}
\end{document}